\documentclass[reqno]{amsart} \usepackage{amscd}
\usepackage{epsf}
\theoremstyle{remark} 
\numberwithin{equation}{section}
\begin{document}
\title[Geometric Quantization]{ Geometric Quantization of the moduli space of the  
vortex equations on a Riemann surface}

\author{Rukmini Dey}

\maketitle

\begin{abstract}
In this note   we quantize the usual symplectic (K\"{a}hler) form  on the vortex moduli space by modifying the Quillen metric of the  Quillen determinant line bundle.
\end{abstract}

\section{The Quantum bundle}

In this note, we highlight a construction done and used in ~\cite{D2}, ~\cite{D3}.

We show that  one can quantize the usual symplectic (in fact, K\"{a}hler) form $\Omega$  on the vortex moduli space, (notation as in ~\cite{D}).   This is done by constructing the line bundle ${\rm det} \bar{\partial}_{A}$   and modifying the Quillen metric by a factor of $e^{- \frac{i}{4 \pi} \int_{M} | \Psi|^2_H \omega} .$ Let us denote this line bundle by ${\mathcal P}.$ This is the procedure followed by Biswas and Raghavendra in constructing determinant bundles on the moduli space of stable triples over a curve, ~\cite{BR}. They used the algebro-geoetric definition of the Quillen bundle.
In ~\cite{D} we contructed the quantum bundle a differently introducing a special section $\Psi_0$, which is a bit unnatural.

We define ${\mathcal P} = {\rm det} (\bar{\partial} + A^{0,1}) =det(\bar{\partial}_A)$ which is well defined on ${\mathcal A} \times \Gamma (L)$ (over every $(A, \Psi)$ the fiber is that of ${\rm det} (\bar{\partial}_A)$). We equip  ${\mathcal P}$ with a modified Quillen metric, namely, we multiply the Quillen metric $e^{-\zeta_A^{\prime}(0)}$ by the factor  $e^{- \frac{i}{4 \pi} \int_{M} | \Psi|^2_H \omega} ,$ where recall $\zeta_A(s)$ is the $zeta$-function 
corresponding to the Laplacian of the $\bar{\partial} + A^{0,1}$ operator. We calculate the curvature for this modified metric on the affine space, ${\mathcal C}$.  
Let $X = ( \alpha_1, \beta)$ and $Y = (\alpha_2, \eta)$ be in $T_{(A, \Psi)} {\mathcal C}$. 
The $zeta$ part of the metric contributes  $-\frac{i}{4 \pi} \int_{M} \alpha_1 \wedge \alpha_2$ to the curvature as is well known. The second part   $e^{-\frac{i}{4 \pi}  \int_{M} | \Psi|^2_H \omega} ,$ contributes to the second part of the curvature form namely, 
$- \frac{i}{4 \pi} \int_{M} (\beta H \bar{\eta} - \bar{\eta} H \beta) \omega$ as follows:

Thus the curvature of ${\mathcal P}$ with the modified Quillen metric is indeed $\frac{i}{ 2 \pi} \Omega$ where 

\begin{eqnarray*}
 \Omega((\alpha_1, \beta), (\alpha_2, \eta)) =  - \frac{1}{2} \int_{M} \alpha_1 \wedge \alpha_2 - \frac{1}{2} \int_{M}( \beta H \bar{\eta} - \bar{\beta} H \eta )  \omega
\end{eqnarray*}
 on the affine space. 

One can show (see for instance ~\cite{D3}), that this line bundle descends to the moduli space of vortices, as long as the descendent of the symplectic form $\Omega$ is integral. The integerality condition is satisfied if the volume of the Riemann surface is integral. This has been showed in ~\cite{MN} when the authors computed the volume of the vortex moduli space.

The quantum bundle is in fact holomorphic with respect to the complex structure on the vortex moduli space and one can take square integrable holomorphic sections as the Hilbert space of the quantization

{\bf Acknowlegement:} I would like to thank Professor Jorgen Andersen, C.T.Q.M., Aarhus University, Denmark, for the very useful discussions. I would like to thank Professor Leon Takhtajan, S.U.N.Y. at Stony Brook, N.Y., for introducing the author to Quillen bundle in the first place.

ICTS-TIFR, Bangalore.
email: rukmini@icts.res.in, rukmini.dey@gmail.com


\begin{thebibliography}{99}
 \bibitem{BF} J.M. Bismut, D.S. Freed: The analysis of elliptic 
families.I. Metrics and connections on determinant bundles; 
Commun. Math. Phys, 106, 159-176 (1986). 
\bibitem{BR} I. Biswas, N. Raghavendra: The determinant bundle on the moduli 
space  of stable triples over a curve; Proc. Indian Acad. Sc. Mat. Sci. 112, 
 no. 3, 367-382, (2002).
\bibitem{D} R. Dey: Geometric prequantization of the moduli space of the 
vortex equations on a Riemann surface; Journal of Mathematical Physics, vol. 47,  issue 10, (2006), page  103501-103508; math-ph/0605025. 
\bibitem{D2} R. Dey: Erratum: Geometric prequantization of the moduli space of the vortex equations 
on a Riemann surface" Journal of Mathematical  
Phys. 50, 119901 (2009)  
\bibitem{D3} R. Dey, V. Mathai, "Holomorphic Quillen determinant bundle on integral compact 
 K\"{a}hler manifolds",  
Quart. J. Math. 64 (2013), 785-794, Quillen Memorial Issue;  arXiv:1202.5213v3
\bibitem{MN} Manton,  N.S., Nasir, S.M., "Volume of vortex moduli spaces",  
Comm. Math. Phys.  199,  no. 3, 591--604 (1998).  
\bibitem{Q} D. Quillen: Determinants of Cauchy-Riemann operators 
over a Riemann surface; Functional Analysis and Its Application, 
19, 31-34 (1985).
\end{thebibliography}
\end{document}